\documentclass[12pt]{article}
\usepackage[dvips]{epsfig}
\usepackage{graphics}
\usepackage{latexsym}
\usepackage{verbatim}
\usepackage{amsmath}
\usepackage{amsthm}
\usepackage{amsfonts}
\usepackage{amssymb}

\newtheorem{theorem}{Theorem}[section]
\newtheorem{lemma}[theorem]{Lemma}
\newtheorem{corollary}[theorem]{Corollary}
\newtheorem{proposition}[theorem]{Proposition}

\def\Remark{\medskip\noindent{\bf Remark: }}
\def\Remarks{\medskip\noindent{\bf Remarks: }}

\newcommand{\vett}{v_{*}}

\newcommand{\ens}[1]{\mathbb{#1}}

\newcommand{\R}{\mathbb{R}}

\def\cst{{\rm cst}} 

\def\var{\varepsilon}
\def\pa{\partial}

\def\signld{\bigskip\bigskip\hspace{80mm}
\vbox{{\sc L. Desvillettes\par\vspace{3mm}
CMLA, ENS Cachan\par
61,  avenue du Président Wilson\par
94235 Cachan Cedex\par
FRANCE\par\vspace{3mm}
e-mail:} desville@cmla.ens-cachan.fr }}

\def\signcm{\bigskip\bigskip\hspace{80mm}
\vbox{{\sc C. Mouhot\par\vspace{3mm}
UMPA, ENS Lyon\par
46 all\'ee d'Italie\par
69364 Lyon Cedex 07\par
FRANCE\par\vspace{3mm}
e-mail:} cmouhot@umpa.ens-lyon.fr }}

\begin{document}

\title{About $L^p$ estimates for the spatially homogeneous Boltzmann equation}

\author{Laurent Desvillettes and Cl\'ement Mouhot}

\hyphenation{bounda-ry rea-so-na-ble be-ha-vior pro-per-ties
cha-rac-te-ris-tic parti-cu-lar}

\date{}

\maketitle

\begin{abstract} 
For the homogeneous Boltzmann equation with (cutoff or non cutoff) 
hard potentials, we prove estimates of 
propagation of $L^p$ norms with a weight $(1 + |x|^2)^{q/2}$ 
($1<p<+\infty$, $q \in \R_+$ large enough), as well as appearance of such weights. 
The proof is based on some new functional inequalities  
for the collision operator, proven by elementary means.
\end{abstract}

\begin{center}
{\bf{ A propos des estimations $L^p$ pour l'\'equation}}
\end{center}
\begin{center}
 {\bf{ de Boltzmann homog\`ene}}
\end{center}
\bigskip

\begin{abstract}
On prouve la propagation de normes $L^p$ avec poids  $(1 + |x|^2)^{q/2}$ et 
l'apparition de tels poids pour l'\'equation de Boltzmann 
homog\`ene dans le cas des potentiels durs (avec ou sans troncature angulaire). 
La d\'emonstration est bas\'ee sur de nouvelles in\'egalit\'es fonctionnelles
pour l'op\'erateur de collision, que l'on prouve par des moyens
\'el\'ementaires
\end{abstract}

\tableofcontents

\section{Introduction}

The spatially homogeneous Boltzmann equation (cf.~\cite{cercignani}) writes
 \begin{equation}\label{eq0}
 \frac{\partial f}{\partial t} (t,v) = Q(f,f)(t,v)\,,
 \end{equation}
where $f(t,\cdot) : \R^{N} \to \R_+$ is the nonnegative
density of particles which at time $t$ move with velocity $v$.
The bilinear
operator in the right-hand side is defined by
 \begin{equation}\label{eq1}
 Q(g,f)(v) = \int_{\R^N} \int_{S^{N-1}} \bigg\{ f(v')\,g(v'_*)
 - f(v)\,g(v_*) \bigg\} \, B \bigg(|v-v_*|, \frac{v-v_*}{|v-v_*|}\cdot \sigma \bigg) \, d\sigma \, dv_* .
 \end{equation}
In this formula, $v',v'_*$ and $v,v_*$
are the velocities of a pair of particles before and after a  collision.
They are defined by
 \begin{equation*}
   v' = \frac{v+\vett}{2} + \frac{|v-\vett|}{2}\,\sigma, \ \ \ 
   \vett' = \frac{v+\vett}{2} - \frac{|v-\vett|}{2}\,\sigma,
 \end{equation*}
where $\sigma\in S^{N-1}$.
\medskip

We concentrate in this work on hard potentials or hard spheres
collision kernels, with or without angular cutoff.
More precisely, we suppose that the collision kernel satisfies the following
\bigskip

\noindent
{\bf{Assumptions}}: The collision kernel $B$ is of the form
 \begin{equation}\label{eq:hyp1}
 B(x,y) = |x|^{\gamma}\, b(|y|), 
 \end{equation}
where 
 \begin{equation} \label{eq:hyp2}
 \gamma\in \, [0,1]
 \end{equation} 
and  
 \begin{equation} \label{eq:hyp3}
 b \in L^{\infty}_{loc}([-1,1[), \ \ \ \ \ 
 b (y) = O_{y\to 1^-} \big((1-y)^{\frac{-(N-2)+\nu}{2}}\big), \ \ \nu>-3.
 \end{equation} 
\bigskip

Note that assumption \eqref{eq:hyp3} is an alternative (and a slighlty less general) formulation to the minimal 
condition necessary for a mathematical treatment of the Boltzmann equation 
identified in~\cite{Vill:new:98,alvi}, namely the requirement
 \begin{equation}\label{eq:hypb} 
 \int_{\ens{S}^{N-1}} b ( \cos \theta ) ( 1 - \cos \theta ) \, d\sigma < +\infty. 
 \end{equation}

\medskip

Then, we wish to consider initial data $f_0 \ge 0$ with
finite mass and energy, such that 
$f_0 (1+|v|^2)^{q/2} \in L^p (\R^N)$ for some $1<p<+\infty$ and 
$q \ge 0$ (notice that entropy is thus automatically finite). 
Existence results under the assumptions of finite mass, energy and entropy  
were obtained in~\cite{ark71} for the case of  hard potentials with cutoff, 
in~\cite{ark:infini:81} for (non cutoff) soft potentials in dimension $3$ under 
the restriction $\gamma \ge -1$, then in~\cite{goudon:grazing:97} 
and~\cite{Vill:new:98} for general kernels (our assumptions 
on the kernel fall in the setting of~\cite{Vill:new:98} for instance). 
Uniqueness however is proved only in the cutoff case (for optimal result see~\cite{mischler-wennberg}) 
and remains an open question in the noncutoff case (except for maxwellian molecules $\gamma=0$, see~\cite{23}).

Propagation of moments in $L^1$ was proven in~\cite{trues} for 
Maxwellian molecules with cutoff. Then, for the case of strictly hard potentials 
with cutoff, it was shown in~\cite{desville-moments} that all polynomial moments 
were created immediately when one of them
of order strictly bigger than $2$ initially existed.
This last restriction was later relaxed in~\cite{Wenn:momt:97}. 

Propagation of moments in $L^p$ was first obtained by
Gustafsson (cf.~\cite{gustafsson1,gustafsson2}) thanks to 
interpolation techniques, under the assumption of
angular cutoff. It was recovered by a simpler and more 
explicit method in~\cite{mv}, thanks to the 
smoothness properties of the gain part of the Boltzmann's 
collision operator discovered by P.-L.~Lions~\cite{lions:regQcutoff:94}. 
As far as appearance of moments in $L^p$ is concerned, the first 
result is due to Wennberg in~\cite{Wenn:L^p:94}, still in the 
framework of angular cutoff. It is precised in~\cite{mv}.
\medskip

In this work, we wish to improve these results  by presenting 
an $L^p$ theory 
 \begin{itemize}
 \item first, which is elementary
 (that is, without abstract interpolations and without using the 
 smoothness properties of Boltzmann's kernel),
 \item secondly, which includes the  non cutoff case,
 \item finally, without assuming too many moments in $L^p$ for the initial datum. 
 \end{itemize}

Our method is reminiscent of recent works by Mischler and 
Rodriguez Ricard \cite{miro} and Escobedo, Lauren\c cot and Mischler \cite{eslami} on the
Smoluchowsky equation. 
\medskip

Let $1<p<+\infty$. We define the weighted $L^p$ space $L^p_q(\R^N)$ by
 \[L^p_q(\R^N) = \bigg\{ f:\R^N \to \R, \quad  ||f||_{L^p_q(\R^N)} <
 +\infty \bigg\}, \]
with its norm
 \[ ||f||_{L^p_q(\R^N)}^p = \int_{\R^N} |f(v)|^p\, \langle v \rangle^{pq}\, dv, \]
and the usual notation $\langle v \rangle = (1+|v|^2)^{1/2}$.
\medskip 

We now state our main theorem
 \begin{theorem}\label{theo:main} 
 Let $B$ be a collision kernel satisfying Assumptions~\eqref{eq:hyp1}, \eqref{eq:hyp2}, 
 \eqref{eq:hyp3} and $q$ such that \\
  (i) $q \in \R_+$ if $\nu > -1$ (integrable angular kernel),\\
  (ii) $pq > 2$ if $\nu \in (-2,-1]$,\\
  (iii) $pq > 4$ if $\nu \in (-3,-2]$,\\
 and $f_0$ be an initial datum in $L^1_{\max(p,2)\,q +2} \cap L^p_q$. 
\smallskip

Then 
  \begin{itemize}
  \item there exists a (weak) solution to the Boltzmann 
  equation~\eqref{eq0} with collision kernel $B$ and initial datum $f_0$
  lying in $L^\infty([0,+\infty);L^p_q(\R^N))$ (with explicit bounds in this space),
  \item if $\gamma >0$, this solution belongs moreover to 
  $L^{\infty}((\tau, +\infty); L^p_r(\R^N))$ for all $\tau >0$ and $r >q$ (still with
  explicit bounds in this space, the blow up near $\tau \sim 0^+$ being at worse polynomial). 
  \end{itemize}
 \end{theorem}

\Remarks We now discuss the assumptions and the conclusion of this theorem. 
\smallskip 

1. Our result cannot hold when the hard potentials are replaced by soft potentials. 
In the case of Maxwellian molecules ($\gamma=0$), we have uniform (in time) bounds 
but no appearance of moments (either in $L^p$ nor 
in $L^1$) occurs. In the case of the so-called ``mollified soft potentials'' with cutoff, 
some bounds growing polynomially in time can be found in~\cite{tovi}, based on the 
regularity property of the gain term of the collision operator.
\smallskip 

2. When the collision kernel $B$ is not a product of a function of $x$ by a function of $y$ 
(as in Assumption~\eqref{eq:hyp1}), it is likely
that Theorem~\ref{theo:main} still holds provided that the behavior of $B$ with respect to
$x$ (when $x \to +\infty$) is that of a nonnegative power and $B$ satisfies 
estimate~\eqref{eq:hyp3} uniformly according to $x$. 
\smallskip

3. The restriction on the weight $q$ is not a technical one which is likely to be discarded 
(at least in our method). Indeed as suggested in~\cite{aldeviwe} the noncutoff collision 
operator behaves roughly like some fractional Laplacian of order $-\nu/2$ and these 
derivatives will in fact be supported by the weight, as we shall see. Notice however 
that there is no condition on $q$ when $\nu > -1$, i.e. in the cutoff case, which is 
coherent with existing results. Note also that the condition $f_0 \in L^1_{2q +2}$
is used only to get the uniformity when $t\to +\infty$ of the estimates. The local (in time)
estimates hold as soon as $f_0 \in L^1_{pq +2}$.
\smallskip 

4. Finally, Theorem~\ref{theo:main} can certainly be improved when the collision kernel in non cutoff. 
In such a case (and under rather not stringent assumption (cf. \cite{aldeviwe})), 
it is possible to show that some smoothness is gained, and some $L^p$ 
regularity will appear even if it does not initially exist. As a consequence, the assumptions of 
Theorem~\ref{theo:main} can certainly be 
somehow relaxed. One can for example compare Theorem~\ref{theo:main} to the results
of~\cite{DVlandau1} for the Landau equation. We also refer to~\cite{desville-wennberg} for ``regularized
hard potentials'' without angular cutoff. 
\medskip

The proof of Theorem~\ref{theo:main} runs as follows. In Section~\ref{sec:estim}, we 
give various bounds for quantities like 
 \[ \int_{\R^N} Q(f,f)(v) \, f^{p-1}(v)\, \langle v \rangle ^{pq}\, dv.\] 
These bounds are applied to the flow of the spatially homogeneous Botzmann 
equation in Section~\ref{sec:flow}, and are sufficient to prove Theorem~\ref{theo:main}, 
except that the bounds may blow up when $t \to +\infty$.
Finally in Section~\ref{sec:unif}, we explain  why
such a blow up cannot take place, and so we conclude the proof of Theorem~\ref{theo:main}.
This last part is the only one which is not self-contained.
It uses an estimate from \cite{mv}.

\section{Functional estimates on the collision operator} \label{sec:estim} 

\begin{figure}[h]
 \epsfysize=4cm
 $$\epsfbox{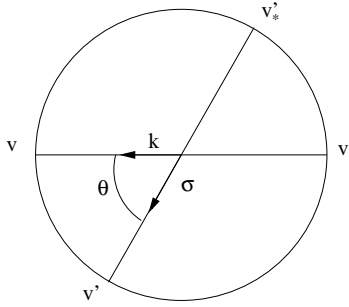}$$
 \caption{Geometry of binary collisions}\label{fig:geom}
 \end{figure}

In the sequel we shall use the parametrization described in figure~\ref{fig:geom}, where 
 \[ \sigma = \frac{v-v_*}{|v-v_*|}, \ \ k = \frac{v' - v' _*}{|v' -v' _*|}, \]
and $\cos \theta = \sigma \cdot k$. 
The range of
$\theta$ is $[0,\pi]$ and $\sigma$ writes
 \[ \sigma = \cos \theta \, k + \sin \theta \, u, \]
where  $u$ belongs to the sphere of $\ens{S}^{N-1}$
orthogonal to $k$ (which is isomorphic to 
$\ens{S}^{N-2}$).
 
Thanks to the change of variable $\theta \mapsto \pi-\theta$ 
which exchanges $v'$ and $v'_*$, the {\em quadratic} collision operator can be written
 \begin{equation*} \nonumber
 Q(f,f)(v) = \int_{\R^N \times \ens{S}^{N-1}} 
 \bigg\{ f(v')\,f(v'_*) - f(v)\,f(v_*) \bigg\} \, B_{sym} (|v-v_*|, \cos\theta) \, d\theta \, dv_*,
 \end{equation*}
where 
 \begin{equation*} 
 B_{sym} (|v-v_*|, \cos \theta) = \big[ B(|v-v_*|, \cos\theta) + 
 B(|v-v_*|, \cos(\pi-\theta)) \big] \, 1_{\cos\theta\geq 0}.
 \end{equation*}
\par
As a consequence, it is enough to consider the case when
$B(|v-v_*|, \cdot)$ has its support included
in $[0,\pi/2]$. 
This is what we shall systematically do in the sequel (Beware that certain propositions are 
written for the bilinear kernel $Q(g,f)$ and not for $Q(f,f)$: they hold only in fact
for the symmetrized collision kernel $B_{sym}$ defined above).
\bigskip

Recalling that
 \[ v' = \frac{v+v_*}2 + \frac{|v-v_*|}2\,\sigma, \]
we use (for all $F$) the formula (cf.~\cite[Section~3, proof of Lemma~1]{aldeviwe})
 \begin{multline}\label{eq:cv}
 \int_{\R^N \times \ens{S}^{N-1}} B(|v-v_*|, \cos\theta) \, F(v') \, dv \, d\sigma \\
 = \int_{\R^N \times \ens{S}^{N-1}} \frac{1}{\cos^N (\theta/2)}
 B\left(\frac{|v-v_*|}{\cos(\theta/2)},\cos\theta \right) \, F(v) \, dv \, d\sigma. 
 \end{multline}

Let us prove a first functional estimate independant on the integrability of 
the angular part of the collision kernel
 \begin{proposition}\label{prop:estim1}
 Let $B$ be a collision kernel satisfying Assumptions~\eqref{eq:hyp1}, \eqref{eq:hyp2}, \eqref{eq:hyp3}. Then, 
 for all $p>1$, $q \in \R$ and $f$ and $g$ nonnegative, we have
  \begin{multline}\label{eq:estim1}
  \int_{\R^N} Q(g,f)(v) \, f^{p-1}(v)\, \langle v \rangle ^{pq} \, dv
  \le \\ \int_{\R^{2N} \times \ens{S}^{N-1}} \, |v-v_*|^{\gamma}
  \, b(\cos\theta) \, 
  \left[ \left(\cos(\theta/2)\right)^{-\frac{N+\gamma}{p'}} - 1 \right] 
  \, \langle v \rangle ^{pq} \, f^p(v) \, g(v_*) \, d\sigma \, dv_* \, dv \\
  +  \int_{\R^{2N} \times \ens{S}^{N-1}}
  \frac1p \, \left(\cos(\theta/2)\right)^{- \frac{N+\gamma}{p'}} \,  |v-v_*|^{\gamma} \,
  b(\cos\theta) \\ \times\,\left[ \langle v' \rangle ^{pq} - \langle v \rangle ^{pq} \right] \, 
  f^p(v) \, g(v_*) \, d\sigma \, dv_* \,dv.
  \end{multline}
 \end{proposition}

\begin{proof}[Proof of Proposition~\ref{prop:estim1}] 
We first observe that thanks to the pre-post collisional 
change of variables (that is, the identity $\int\int\int F(v,v_*,\sigma) \, d\sigma \, dv_* \, dv =
\int\int\int F(v',v'_*,\sigma) \, d\sigma \, dv_* \, dv $):
 \begin{eqnarray*}
 && \int_{\R^N} Q(g,f)(v) \, f^{p-1}(v)\, \langle v \rangle ^{pq} \, dv \\ 
 &=& \int_{\R^{2N} \times \ens{S}^{N-1}} \Big\{  g(v'_*) \, f(v')
 - \, g(v_*)\,f(v) \Big\} \, f^{p-1}(v) \, \langle v \rangle ^{pq} \, |v-v_*|^{\gamma}
 \, b(\cos\theta)\, d\sigma \, dv_* \, dv \\
 &=&  \int_{\R^{2N} \times \ens{S}^{N-1}} \Big[\langle v' \rangle ^{pq} \,f^{p-1}(v') \,
 f(v) \, g(v_*) \\
 && \ \ \ \ \ \ \ \ \ \ \ \ \ \ \ \ \ \ \ \ \ \ \ \ \ - \, \langle v \rangle ^{pq} \, f^p(v) \, g(v_*) \Big]
 |v-v_*|^{\gamma} \, b(\cos\theta)\, d\sigma \, dv_* \, dv.
 \end{eqnarray*}
According to Young's inequality, for all $\mu \equiv \mu(\theta) > 0$, 
 \begin{equation*} 
 f^{p-1}(v')\,f(v) = \bigg( \frac{f(v')}{\mu^{1/p}} \bigg)^{p-1}\, 
 (\mu^{1 - 1/p}\, f(v))
 \le \left(1 - \frac1p\right) \, \mu^{-1}\, f^p(v') + \frac1p\, \mu^{p-1}\,f^p(v) ,
 \end{equation*}
so that
 \begin{multline*} 
 \int_{\R^N} Q(g,f)(v) \, f^{p-1}(v)\, \langle v \rangle ^{pq} \, dv  
 \le \int_{\R^{2N} \times \ens{S}^{N-1}} \bigg[ \left(1 - \frac1p\right) \, 
 \mu^{-1}\, \langle v' \rangle ^{pq} \, f^p(v') \\
 + \, \frac1p\, \mu^{p-1}\, \langle v' \rangle ^{pq} \, f^p(v) - \langle v \rangle ^{pq} \, 
 f^p(v) \bigg] \, g(v_*)\, |v-v_*|^{\gamma}
 \, b(\cos\theta) \, d\sigma \, dv_* \, dv.
 \end{multline*}
We now use (for a given $v_*$, $\theta$) formula~\eqref{eq:cv} for the first
term in this integral. We get
 \begin{eqnarray*}
 && \int_{\R^N} Q(g,f)(v) \, f^{p-1}(v)\, \langle v \rangle ^{pq} \, dv 
 \le \int_{\R^{2N} \times \ens{S}^{N-1}} 
 \bigg[ \left(1 - \frac1p\right) \, \mu^{-1}\, \langle v \rangle ^{pq} 
 \, (\cos(\theta/2))^{-N-\gamma}\, f^p(v) \\
 && \ \ \ \ \ \ \ \ \ \ \ \ \ \ 
 + \, \frac1p\, \mu^{p-1}\, \langle v' \rangle ^{pq} \, f^p(v) - \langle v \rangle ^{pq} \, 
 f^p(v) \bigg] \, g(v_*) \, |v-v_*|^{\gamma} \, b(\cos\theta) \, d\sigma \, dv_* \, dv \\ 
 && \ \ = \int_{\R^{2N} \times \ens{S}^{N-1}} \langle v \rangle ^{pq} \, |v-v_*|^{\gamma}
 \, b(\cos\theta) \,  f^p(v)\,  g(v_*) \\
 && \ \ \ \ \ \ \ \ \ \ \ \ \ \ \ \ \ 
     \times \, \bigg[ \left(1 - \frac1p\right) \, \mu^{-1}\, (\cos(\theta/2))^{-N-\gamma}
 + \frac1p\, \mu^{p-1}  - 1 \bigg] \, d\sigma \, dv_* \, dv \\ 
 && \ \ \ \ \ + \int_{\R^{2N} \times \ens{S}^{N-1}} \frac1p\, \mu^{p-1} \, 
 |v-v_*|^{\gamma} \, b(\cos\theta) \, f^p(v) \, g(v_*)
 \times \, \left[ \langle v' \rangle ^{pq}
 - \langle v \rangle ^{pq} \right] \,  d\sigma \, dv_* \, dv.
 \end{eqnarray*}
We now take the optimal  $\mu = \mu(\theta) > 0$. This amounts to consider
 \[ \mu(\theta) = (\cos(\theta/2))^{- \frac{N+\gamma}{p}}.\]
 In this way, we get estimate 
(\ref{eq:estim1}). 
\end{proof}

\Remark With the same idea, one could easily obtain 
  \begin{multline*}
 \int_{\R^N} Q(g,f)(v) \, f^{p-1}(v)\, \langle v \rangle ^{pq} \, dv \\
 =  \int_{\R^{2N} \times \ens{S}^{N-1}} \langle v \rangle ^{pq} \, |v-v_*|^{\gamma}
 \, b(\cos\theta) \,  f^p(v)\,  g(v_*)  \\
 \times\,\bigg[ \left(1 - \frac1p\right) \, \mu^{-1}\, (\cos(\theta/2))^{-N-\gamma}
 + \frac1p\, \mu^{p-1}\, (\cos(\theta/2))^{pq} - 1 \bigg] \, d\sigma \, dv_* \, dv \\ 
 + \, \int_{\R^{2N} \times \ens{S}^{N-1}} \frac1p\, \mu^{p-1} \, 
 |v-v_*|^{\gamma} \, b(\cos\theta) \, f^p(v) \, g(v_*) 
 \, \big[ \langle v' \rangle ^{pq}
 - (\cos(\theta/2)) ^{pq} \langle v \rangle ^{pq} \big] \,  d\sigma \, dv_* \, dv,
 \end{multline*}
so that taking the optimal $\mu$ given by
 \[ \mu(\theta) = (\cos(\theta/2))^{- \frac{N+\gamma}{p} -q}, \]
the following inequality holds:
  \begin{multline*}
  \int_{\R^N} Q(g,f)(v) \, f^{p-1}(v)\, \langle v \rangle ^{pq} \, dv \le \\ 
  \int_{\R^{2N} \times \ens{S}^{N-1}} \langle v \rangle ^{pq} \, |v-v_*|^{\gamma}
  \, b(\cos\theta) \, 
  \left[ \left(\cos(\theta/2)\right)^{q-\frac{N+\gamma}{p'}} - 1 \right] 
  \, f^p(v) \, g(v_*) \, d\sigma \, dv_* \, dv \\
  +  \int_{\R^{2N} \times \ens{S}^{N-1}}
  \frac1p \, \left(\cos(\theta/2)\right)^{-q(p-1) - \frac{N+\gamma}{p'}} \,  |v-v_*|^{\gamma} \,
  b(\cos\theta) \\ 
  \big[ \langle v' \rangle ^{pq} - 
  (\cos(\theta/2)) ^{pq} \langle v \rangle ^{pq} \big] \, 
  f^p(v) \, g(v_*) \, d\sigma \, dv_* \,dv.
  \end{multline*}
If $q$ is big enough, i.e. such that 
 \begin{equation}\label{eq:badhyp}
 q-\frac{N+\gamma}{p'} >0,
 \end{equation} 
the first term is strictly negative, and 
some estimates (in the same spirit as in Lemma~\ref{lem:sym} below)
on the term $\big[ \langle v' \rangle ^{pq} - (\cos(\theta/2)) ^{pq} \langle v \rangle ^{pq} \big]$ 
for small and large angles $\theta$ would yield {\em directly} to 
 \begin{eqnarray*} 
 \int_{\R^N} Q(g,f)(v) \, f^{p-1}(v)\, \langle v \rangle ^{pq} \, dv \le 
 &-& \, C \, \int_{\R^N} g (v_*)\,dv_* \,  \int_{\R^N} f^{p}(v) \, \langle v \rangle ^{pq + \gamma} \, dv \\
 &+& \, D\, \int_{\R^N} g(v_*) \,  \langle v \rangle^{pq + \gamma} \, dv_*  
 \, \int_{\R^N} f^{p}(v) \,dv \\
 &+& \, D\, \int_{\R^N} g(v_*) \,  \langle v \rangle ^2 \, dv_* 
 \, \int_{\R^N} f^{p}(v) \, \langle v \rangle ^{pq} \,dv .
 \end{eqnarray*}
We do not follow in the sequel this line of ideas because we don't want
to assume~\eqref{eq:badhyp}.
We rather choose to make a global splitting between the small and large angles $\theta$.
\medskip

We now deduce from Proposition~\ref{prop:estim1} a corollary enabling to bound
 \[ \int_{\R^N} Q(g,f)(v) \, f^{p-1}(v)\, \langle v \rangle ^{pq} \, dv \] 
in terms of weighted $L^1$ and $L^p$ norms of $f$ and $g$. 
Note that this corollary is almost obvious to prove
when the collision kernel is integrable (cutoff case).
 
\begin{corollary}\label{coro:estim2}
 Let $B$ be a collision kernel satisfying Assumptions~\eqref{eq:hyp1}, \eqref{eq:hyp2}, 
 \eqref{eq:hyp3}. We consider
  $f$ and $g$ nonnegative and $q \in \R$. We suppose moreover that
 $pq \ge 2$ if $\nu \in (-2,-1]$ and $pq \ge 4$ if $\nu \in (-3,-2]$. Then, 
  \begin{equation} \label{eq:fonc}
  \int_{\R^N} Q(g,f)(v) \, f^{p-1}(v)\, \langle v \rangle ^{pq} \, dv \le 
  C _{p,N,\gamma} (b) \, \|g\|_{L^1 _{pq + \gamma}} \, \|f\|^p _{L^p _{q+\gamma/p}}  
  \end{equation}
 where 
  \[ C _{p,N,\gamma} (b) = 
 \mbox{cst} \, (p,N,\gamma) \, 
 \left( \int_{\ens{S}^{N-1}} b(\cos \theta) \, ( 1- \cos\theta ) \,  d\sigma \right), \]
and $\mbox{cst} \, (p,N,\gamma)$ is a computable constant depending on $p$, $N$ and $\gamma$.
 \end{corollary}

\Remark Since the non cutoff collision operator behaves roughly like some fractional 
Laplacian of order $-\nu/2$, one could wonder how  
a functional inequality which does not contain derivatives of the function $f$ can hold. 
The answer is that the pre-post collisional change of variable and formula~\eqref{eq:cv} 
(which play here the role played by integration by part for differential operators) 
 allow to transfer the derivatives on the weight function $\langle v \rangle ^{pq}$. 
This also explains why the restriction on the weight exponent $q$ depends on the order 
$\nu$ of the angular singularity.
\medskip

\begin{proof}[Proof of Corollary~\ref{coro:estim2}]
Estimate~\eqref{eq:estim1} can be written 
 \begin{equation*}
 \int_{\R^N} Q(g,f)(v) \, f^{p-1}(v)\, \langle v \rangle ^{pq} \, dv 
 \le I_1 + I_2 + I_3 \,, 
 \end{equation*} 
where 
 \begin{eqnarray*} 
 I_1 &=&  \int_{\R^{2N} \times \ens{S}^{N-1}}  |v-v_*|^{\gamma}
  \, b(\cos\theta) \, 
  \left[ \left(\cos(\theta/2)\right)^{-\frac{N+\gamma}{p'}} - 1 \right] 
  \, \langle v \rangle ^{pq} \, f^p(v) \, g(v_*) \, d\sigma \, dv_* \, dv, \\
 I_2 &=&  \int_{\R^{2N} \times \ens{S}^{N-1}} 
  \frac1p \, \left[ \left(\cos(\theta/2)\right)^{- \frac{N+\gamma}{p'}} -1 \right] \,  |v-v_*|^{\gamma} \,
  b(\cos\theta) \,\left[ \langle v' \rangle ^{pq} - \langle v \rangle ^{pq} \right] \, 
  f^p(v) \, g(v_*) \, d\sigma \, dv_* \,dv, \\ 
 I_3 &=&  \int_{\R^{2N} \times \ens{S}^{N-1}} 
  \frac1p \,  |v-v_*|^{\gamma} \,
  b(\cos\theta) \,\left[ \langle v' \rangle ^{pq} - \langle v \rangle ^{pq} \right] \, 
  f^p(v) \, g(v_*) \, d\sigma \, dv_* \,dv. 
 \end{eqnarray*} 
Then the two first terms are easily estimated thanks to the formula
 \[ \left[ \left(\cos(\theta/2)\right)^{-\frac{N+\gamma}{p'}} - 1 \right] 
    \sim_{\theta \to 0} \frac{N+\gamma}{4 p'} (1-\cos \theta ). \]
For the last one, we shall need the following lemma, which takes advantage of the 
symmetry properties of the collision operator: 
 \begin{lemma}\label{lem:sym}
 For all $\alpha \ge 1$, 
  \begin{equation} 
  \bigg| \int_{u\in  \ens{S}^{N-2}} 
  \big[ \langle v' \rangle^{2 \alpha} - \langle v \rangle ^{2 \alpha} \big] \, du \bigg|
  \le C_{\alpha}\, (\sin \theta/2) \, \langle v \rangle ^{2\alpha} \, 
  \langle v_* \rangle ^{2 \alpha},
  \end{equation}
 and for all $\alpha \ge 2$, 
  \begin{equation} 
  \bigg| \int_{u\in  \ens{S}^{N-2}} 
  \big[ \langle v' \rangle^{2 \alpha} - \langle v \rangle ^{2 \alpha} \big] \, du \bigg|
  \le C_{\alpha}\, (\sin \theta/2)^2\, \langle v \rangle ^{2 \alpha } \, 
  \langle v_* \rangle ^{2 \alpha }.
  \end{equation}
 \end{lemma}
\Remark This lemma is reminiscent of the symmetry properties 
used in the ``cancellation lemma'' in~\cite{alvi} and~\cite{aldeviwe} in order 
to give sense to the Boltzmann collision operator for strong angular singularity 
(i.e. $\nu \in (-3,-2]$).  
\begin{proof}[Proof of Lemma~\ref{lem:sym}]
We note that since
 \[ |v'|^2 = |v|^2\, \cos^2\theta/2 + |v_*|^2 \, \sin^2\theta/2 + 
 2 \cos\theta/2 \, \sin\theta/2 \, |v-v_*|\, u\cdot v_* ,\]
if one introduces (for $x \in [0, \sqrt{2}/2]$) the function 
 \[ R_\alpha (x) = \int_{u\in  \ens{S}^{N-2}} \left[ \Big(1 + |v|^2 \, (1 - x^2) + |v_*|^2 \, x^2 
 + 2 x\,\sqrt{1-x^2} \, |v-v_*|\, u\cdot v_* \Big)^{\alpha}  - (1+
 |v|^2)^{\alpha} \right] \, du, \]
we get 
 \[ \int_{u\in  \ens{S}^{N-2}} \left[ (1+|v'|^2)^{\alpha} - 
 (1+|v|^2)^{\alpha} \right] \, du = R_\alpha (\sin\theta/2). \]
But thanks to the change of variables $u \to -u$, we see that $R_\alpha$ is even.
Noticing also that $R_\alpha(0) = 0$, we use the identities 
 \begin{eqnarray*} 
 R_\alpha(x) &=& x \, \int_{0}^1 \, R' _\alpha(s\,x) \, ds, \\
 R_\alpha(x) &=& x^2 \, \int_{0}^1 (1-s) \, R'' _\alpha(s\,x) \, ds.
 \end{eqnarray*}
We compute
 \begin{eqnarray*}
 R' _\alpha(x) &=& \alpha \, \int_{u\in  \ens{S}^{N-2}} \Big(-2x\,|v|^2 + 2x\,|v_*|^2 
 + 2 \, (1-x^2)^{1/2}\, |v-v_*|\, u\cdot v_* \\ 
 && - 2 x^2\, 
 (1-x^2)^{-1/2} \, |v-v_*|\, u\cdot v_* \Big) \\
 && \times \, \Big(1 + |v|^2\, (1 - x^2) + |v_*|^2 \, x^2 
 + 2 x\,\sqrt{1-x^2} \, |v-v_*|\, u\cdot v_* \Big)^{\alpha - 1} \, du 
 \end{eqnarray*}
and 
 \begin{eqnarray*}
 R'' _\alpha(x) &=&  \alpha \, (\alpha - 1) \, \int_{u\in  \ens{S}^{N-2}}  
 \Big( -2x\,|v|^2 + 2x\,|v_*|^2 + 2 \,(1-x^2)^{1/2}\, |v-v_*|\, u\cdot v_* \\ 
 && - 2 x^2\,
 (1-x^2)^{-1/2} \, |v-v_*|\, u\cdot v_* \Big)^2 \\ 
 && \times \, \Big(1 + |v|^2\, (1 - x^2) + |v_*|^2 \, x^2 
 + 2 x\,\sqrt{1-x^2} \, |v-v_*|\, u\cdot v_*\Big)^{\alpha - 2} \, du \\ 
 && +  \alpha \, \int_{u\in  \ens{S}^{N-2}}  \Big(-2|v|^2 + 2\,|v_*|^2 
 - 2 x \,(1-x^2)^{-1/2}  \\ 
 && |v-v_*|\, u\cdot v_* 
 - 2\, |v-v_*|\, u\cdot v_* \, \big( 2x\,(1-x^2)^{-1/2} + x^3\,(1-x^2)^{-3/2} \big) \Big) \\
 && \times \, \Big(1 + |v|^2\, (1 - x^2) + |v_*|^2 \, x^2 
 + 2 x\,\sqrt{1-x^2} \, |v-v_*|\, u\cdot v_* \Big)^{\alpha - 1} \, du.
 \end{eqnarray*}
Then, for $x \in [0, \sqrt{2}/2]$, if $\alpha \ge 1$, we get  
 \begin{equation*}
 | R' _\alpha(x)| \le  C_{\alpha} \, \langle v \rangle^{2\alpha} \, \langle v_* \rangle^{2\alpha}, 
 \end{equation*}
and if $\alpha \ge 2$,  
 \begin{equation*}
 | R'' _\alpha(x)| \le  C_{\alpha} \, \langle v \rangle^{2\alpha} \, \langle v_* \rangle^{2\alpha}.
 \end{equation*}
This concludes the proof of Lemma~\ref{lem:sym}.
\end{proof}
Let us come back to the proof of Corollary~\ref{coro:estim2}. We have 
 \[ I_3 =  \int_{\R^{2N}} \int_0 ^\pi 
  \frac1p \,  |v-v_*|^{\gamma} \, b(\cos \theta) \, R_\alpha(\sin\theta/2) \, 
  (\sin \theta)^{N-2} \, f^p(v) \, g(v_*) \, d\theta \, dv_* \,dv \] 
for $\alpha = (pq)/2$. Lemma~\ref{lem:sym} and the equality
 \[ (\sin \theta/2)^2 = \frac{(1-\cos \theta)}{2} \]
conclude the proof.
\end{proof}

We now turn to an estimate which holds when the (angular part of the) 
collision kernel has its support  in $[\theta_0,\pi/2]$ for some $\theta_0>0$. As we shall see
later on, this term is the ``dominant part'' of the same quantity when the (angular part of 
the) collision kernel has its support in $[0,\pi/2]$.
 
 \begin{proposition}\label{prop:estim3}
 Let $B$ satisfy Assumptions~\eqref{eq:hyp1}, \eqref{eq:hyp2}, \eqref{eq:hyp3}. We suppose 
 moreover that $b$ has its support in $[\theta_0,\pi/2]$. Then, 
 for all $p>1$, $q \ge 0$ and $f$ nonnegative with bounded 
 $L^1 _{pq+2}$ norm, we have
  \begin{equation}
  \int_{\R^N} Q(f,f)(v) \, f^{p-1}(v)\, \langle v \rangle ^{pq} \, dv \le 
  C^+ (b) \, \|f\|^p _{L^p _q} - K^- (b) \, \|f\|^p _{L^p _{q+\gamma/p}}
  \end{equation}
 with 
  \[ C^+ (b) = C^+ \, \left( \int_{\ens{S}^{N-1}} b \, d\sigma \right), \ \ \ \ 
     K^- (b) = K^- \, \left( \int_{\ens{S}^{N-1}} b \, d\sigma \right), \]
 where $C^+$, $K^-$ are strictly positive constants. Both depend on an upper bound on 
 $\|f\|_{L^1 _{pq+2}}$ and on a lower bound on $\|f\|_{L^1}$; 
 $C^+$ also depends on $\theta_0$. 
 \end{proposition}
\Remark This estimate could be deduced from the results of~\cite{mv}, but we shall give here an 
elementary self-contained proof, in the same spirit as that of the proof of Proposition \ref{prop:estim1}.
\medskip

\begin{proof}[Proof of Proposition~\ref{prop:estim3}]
Let us write the quantity to be estimated 
 \begin{eqnarray*}
 \int_{\R^N} Q(f,f)(v) \, f^{p-1}(v)\, \langle v \rangle ^{pq} \, dv 
 &\le& \int_{\R^N} Q^+(f,f)(v) \, f^{p-1}(v)\, \langle v \rangle ^{pq} \, dv \\
 &&   - \int_{\R^N} Q^-(f,f)(v) \, f^{p-1}(v)\, \langle v \rangle ^{pq} \, dv, 
 \end{eqnarray*}
splitting as usual the operator between its gain and loss parts (remember that the small angles 
have been cutoff). On one hand, using $|v-v_*|^\gamma \ge [ \langle v\rangle^\gamma
- \cst \, \langle v_* \rangle^\gamma]$ we get 
 \[ - \int_{\R^N} Q^-(f,f)(v) \, f^{p-1}(v)\, \langle v \rangle ^{pq} \, dv 
    \le -K_0 \, \|b\|_{L^1(\ens{S}^{N-1})} \, \|f\|^p _{L^p _{q+\gamma/p}} 
        + C_0 \, \|b\|_{L^1(\ens{S}^{N-1})} \, \|f\|^p _{L^p _q} \] 
for some constant $K_0 >0$ depending on a lower bound on $\|f\|_{L^1}$ 
and $C_0 >0$ depending on an upper bound on the $\|f\|_{L^1 _\gamma}$.  
On the other hand, 
 \[ \int_{\R^N} Q^+(f,f)(v) \, f^{p-1}(v)\, \langle v \rangle ^{pq} \, dv 
    = \int_{\R^{2N} \times \ens{S}^{N-1}} f'_* \, f' \, f^{p-1} 
    \, \langle v \rangle ^{pq} \, B \, dv \, dv_* \, d\sigma \]
can be split into 
 \begin{eqnarray*} 
 I_1 &=& \int_{\R^{2N} \times \ens{S}^{N-1}} f' _* \, (f\, j_r)'  \, f^{p-1}
    \, \langle v \rangle ^{pq} \, B \, dv \, dv_* \, d\sigma, \\
 I_2 &=& \int_{\R^{2N} \times \ens{S}^{N-1}} f' _* \, (f\, j_{r^c})' \, f^{p-1}
    \, \langle v \rangle ^{pq} \, B \, dv \, dv_* \, d\sigma,
 \end{eqnarray*}
with $j_r (v) = 1_{|v| \le r}$ and $j_{r^c} = 1 -j_r$. This means that we treat separately  
large and small velocities. Then 
 \begin{eqnarray*}  
 I_1 &=& \int_{\R^{2N} \times \ens{S}^{N-1}} f_* \, (f\, j_r) \, (f') ^{p-1}  
    \, \langle v' \rangle ^{pq} \, B \, dv \, dv_* \, d\sigma \\
    &\le& \int_{\R^{2N} \times \ens{S}^{N-1}} f_* \, 
    \left[ \left(1 - \frac1p\right) \, \mu_1 ^{-1}\, f^p(v') + \frac1p \, \mu_1 ^{p-1}\,(f \, j_r)^p(v) \right] 
    \, \langle v' \rangle ^{pq} \, B \, dv \, dv_* \, d\sigma \\
    &\le& \|b\|_{L^1(\ens{S}^{N-1})} \, \Bigg[
    \left(1 - \frac1p\right) \, \mu_1 ^{-1}\, (\cos \pi/4)^{-N-\gamma} \, 
          \|f\|_{L^1 _\gamma} \, \|f\|^p _{L^p _{q +\gamma/p}} \\
    && \ \ \ \ \ \ \ \ \ \ \ \ \ \ \ \ \ \ + \frac1p \, \mu_1 ^{p-1}\, \|f\|_{L^1 _{pq + \gamma}} 
          \, \|f \, j_r\|^p _{L^p _{q +\gamma/p}} \Bigg], 
 \end{eqnarray*}        
and thus 
 \begin{multline} \label{eq:i1}
 I_1 \ \le \ \|b\|_{L^1(\ens{S}^{N-1})} \, \Bigg[ \left(1 - \frac1p\right) \, 
                     \mu_1 ^{-1}\, (\cos \pi/4)^{-N-\gamma} \, 
         \|f\|_{L^1 _\gamma} \, \|f\|^p _{L^p _{q +\gamma/p}} \\ 
         + \frac1p \, \mu_1 ^{p-1} \, r^{\gamma} \, 
           \|f\|_{L^1 _{pq + \gamma}} \, \|f\|^p _{L^p _q}\Bigg]. 
 \end{multline} 
As for $I_2$, we get 
 \begin{equation*}  
 I_2 = \int_{\R^{2N} \times \ens{S}^{N-1}} f' \, (f\, j_{r^c})' _* \, f^{p-1}
    \, \langle v \rangle ^{pq} \, \tilde{B} \, dv \, dv_* \, d\sigma
 \end{equation*}
thanks to the change of variable $\sigma \to -\sigma$. Now $\tilde{B}$ has compact 
support in $[\pi/2,\pi-\theta_0]$. Then we compute
 \begin{eqnarray*}  
 I_2 &=& \int_{\R^{2N} \times \ens{S}^{N-1}} (f\, j_{r^c})_* \, f \, (f') ^{p-1}  
    \, \langle v' \rangle ^{pq} \, \tilde{B} \, dv \, dv_* \, d\sigma \\
    &\le& \int_{\R^{2N} \times \ens{S}^{N-1}} (f\, j_{r^c})_* \, 
    \left[ \left(1 - \frac1p\right) \, \mu_2 ^{-1}\, f^p(v') + \frac1p \, \mu_2 ^{p-1}\,f^p(v) \right] 
    \, \langle v' \rangle ^{pq} \, \tilde{B} \, dv \, dv_* \, d\sigma \\
    &\le& \|b\|_{L^1(\ens{S}^{N-1})} \, \Bigg[ 
          \left(1 - \frac1p\right) \, \mu_2 ^{-1}\, (\sin \theta_0/2)^{-N-\gamma} \,  
          \|f \, j_{r^c}\|_{L^1 _\gamma} \, \|f\|^p _{L^p _{q +\gamma/p}} \\
    && \ \ \ \ \ \ \ \ \ \ \ \ \ \ \ \ \ \ \ \ \ \ \ \ \ \ \ \ \ \ \ \ \ \ \ \ \ \ \ \ 
          + \frac1p \, \mu_2 ^{p-1}\, \|f \, j_{r^c}\|_{L^1 _{pq + \gamma}} 
          \, \|f\|^p _{L^p _{q +\gamma/p}} \Bigg] 
 \end{eqnarray*}  
by using again formula~\eqref{eq:cv} and thus 
 \begin{eqnarray} \label{eq:i2}\nonumber 
 I_2 &\le& \|b\|_{L^1(\ens{S}^{N-1})} \, \Big[ 
         \left(1 - \frac1p\right) \, \mu_2 ^{-1} \, (\sin \theta_0/2)^{-N-\gamma} \, 
         (1+r^2)^{(\gamma-2)/2} \, \|f\|_{L^1 _2} \, \|f\|^p _{L^p _{q +\gamma/p}} \\  
     && \ \ \ \ \ \ \ \ \ \ \ \ \ \ \ \ \ \ \ \ \ \ \ \ \ \ \ \ \ 
         + \frac1p \, \mu_2 ^{p-1} \,  
         \|f\|_{L^1 _{pq + \gamma}} \, \|f\|^p _{L^p _{q+\gamma/p}} \Big].  
 \end{eqnarray} 
Gathering~\eqref{eq:i1} and~\eqref{eq:i2}, we obtain for the gain part
 \begin{multline*}
 \int_{\R^N} Q^+(f,f)(v) \, f^{p-1}(v)\, \langle v \rangle ^{pq} \, dv \le \\
     \|b\|_{L^1(\ens{S}^{N-1})} \,\left[ \frac1p \, 
     \mu_1 ^{p-1} \, (1+r^2)^{\gamma/2} \, \|f\|_{L^1 _{pq + \gamma}} \right] \|f\|^p _{L^p _q} \\
     + \|b\|_{L^1(\ens{S}^{N-1})} \, \Bigg[ \left(1 - \frac1p\right) 
         \, \mu_1 ^{-1}\, (\cos \pi/4)^{-N-\gamma} \\
     + \left(1 - \frac1p\right) \, \mu_2 ^{-1}\, (\sin \theta_0/2)^{-N-\gamma} \, (1+r^2)^{(\gamma-2)/2} 
     + \frac1p \, \mu_2 ^{p-1} \Bigg] \, \|f\|_{L^1 _{pq+ \gamma}} \, \|f\|^p _{L^p _{q +\gamma/p}}. 
 \end{multline*}
For some $\theta_0>0$ fixed, one can first  choose $\mu_2$ small enough, then $r$ big enough (remember 
that $\gamma-2 <0$), then $\mu_1$ big enough, in such a way that
 \begin{equation*}
 \left[ \left(1 - \frac1p\right) \, \mu_1 ^{-1}\, (\cos \pi/4)^{-N-\gamma} 
       + \left(1 - \frac1p\right) \, \mu_2 ^{-1}\, (\sin \theta_0/2)^{-N-\gamma} \, r^{\gamma-2} 
       + \frac1p \, \mu_2 ^{p-1} \right]\, \|f\|_{L^1 _{pq+ \gamma}} \le \frac{K_0}{2}. 
 \end{equation*}
We thus get the wanted estimate by combining the estimates for the gain part and the loss part.
\end{proof}
 
We now can gather Corollary~\ref{coro:estim2} 
with Proposition~\ref{prop:estim3} in order to get the
 \begin{proposition}\label{prop:estim5} 
 Let $B$ satisfy Assumptions~\eqref{eq:hyp1}, \eqref{eq:hyp2}, \eqref{eq:hyp3},
 $p$ belong to $(1,+\infty)$, and $q\ge 0$.  We suppose moreover that
 $pq \ge 2$ if $\nu \in (-2,-1]$ and $pq \ge 4$ 
 if $\nu \in (-3,-2]$. Then, for $f$ nonnegative with bounded $L^1 _{pq+2}$ norm, we have
  \begin{equation}
  \int_{\R^N} Q(f,f)(v) \, f^{p-1}(v) \, \langle v \rangle ^{pq} \, dv \le 
  C^+ \, \|f\|^p _{L^p _q} 
  - K^- \, \|f\|^p _{L^p _{q+\gamma/p}}
  \end{equation}
 for some positive constants $C^+$ and $K^-$, depending on an upper bound on
 $\|f\|_{L^1 _{pq+2}}$ and on a lower bound on $\|f\|_{L^1}$.
 \end{proposition}

\begin{proof}[Proof of Proposition~\ref{prop:estim5}]
The proof is straightforward and based on a splitting of $b$ of the form
 \begin{equation} 
 b = b_c ^{\theta_0} + b_r ^{\theta_0},
 \end{equation}
where $b_c ^{\theta_0}  = b \, 1_{\theta \in [\theta_0,\pi/2]}$ stands for the ``cutoff'' part, 
$b_r ^{\theta_0}= 1-b_c ^{\theta_0}$ for the remaining part, and $\theta_0 \in (0,\pi/2]$ is 
some fixed positive angle.
We split the corresponding collision operator as  $Q = Q_c + Q_r$. It remains then to 
apply Corollary~\ref{coro:estim2} to
 \[ \int_{\R^N} Q_r (f,f)(v) \, f^{p-1}(v) \, \langle v \rangle ^{pq} \, dv \]
and Proposition~\ref{prop:estim3} to
 \[ \int_{\R^N} Q_c (f,f)(v) \, f^{p-1}(v) \, \langle v \rangle ^{pq} \, dv. \]
Observing that
 \[ \int_{\ens{S}^{N-1}} b_r ^{\theta_0} (\cos \theta) \, ( 1- \cos\theta ) \,  d\sigma 
    \to_{\theta_0 \to 0} 0, \]
we see that the term corresponding to  $Q_r$ can be absorbed 
 by the damping (nonpositive) part
of  $Q_c$, for $\theta_0$ small enough. 
\end{proof}

\section{Application to the flow of the equation}\label{sec:flow}

In this section, we denote by $K$ any strictly 
positive constant which can be replaced by a smaller
strictly positive constant, and by $C$ any constant which can be
replaced by a larger constant. We precise the dependance with respect 
to time when this is useful.
\medskip

We now prove Theorem~\ref{theo:main} without trying to get bounds which are uniform when $t\to +\infty$. 
We notice that  a solution $f(t, \cdot)$ at time $t \ge 0$ of the Boltzmann equation (given by
the results of \cite{ark71},~\cite{ark:infini:81} and~\cite{Vill:new:98}) 
satisfies:
 \[ \frac{d}{dt}  \int_{\R^N}  f^p(v) \, \langle v \rangle^{pq} \, dv = 
 p \, \int_{\R^N} Q(f,f)(v) \, f^{p-1}(v)\, \langle v \rangle^{pq} \, dv.\]
We also recall that (under our assumptions on the initial datum), such a solution $f(t, \cdot)$ 
has a constant mass $||f(t, \cdot)||_{L^1}$. The $L^p_q$ integrability of 
the initial datum $f_0$ implies that this initial datum has bounded 
entropy, then the $H$-Theorem ensures that the entropy remains uniformly bounded 
for all times (by the initial entropy). Also its moment of order $2+pq$ in $L^1$ is propagated 
and remains uniformly bounded for all times with explicit constant (see for 
instance~\cite{Wenn:momt:97}).

\par
  Then Propostion~\ref{prop:estim5} gives the following a priori differential inequality, 
 \begin{equation}\label{eq:apriori}
 \frac{d}{dt} \, \| f\|^p _{L^p_q} \le 
 C \, \|f\|^p _{L^p _q} 
  - K \, \|f\|^p _{L^p _{q+\gamma/p}}.
 \end{equation}
 In particular, 
\begin{equation}\label{ghgh}
 \frac{d}{dt} \, \| f\|_{L^p_q}^p \le C \, \| f\|_{L^p_q}^p.
\end{equation}
According to Gronwall's lemma, the norm $ \| f\|_{L^p_q}$ remains bounded 
(on all intervals $[0,T]$ for $T>0$) if it is initially finite. 

Let us now turn to the question of appearance of higher moments in $L^p$ (when $\gamma >0$). 
Let $r>0$. Using H\"older's inequality, we see that 
 \[ \|f\|_{L^p _r} \le \|f\|^\theta _{L^p _{q_1}} \, \|f\|^{1-\theta} _{L^p _{q_2}} \]
with $r = \theta q_1 + (1-\theta) q_2$. Thus with $q_2 =0$ and $q_1 = r +\gamma/p$, 
we get 
 \[ \|f\|_{L^p _r} \le \|f\|^{\frac{r}{r+\gamma/p}} _{L^p _{r+\gamma/p}} 
                       \, \|f\|^{\frac{\gamma/p}{r+\gamma/p}} _{L^p} \]
Therefore, 
 \[   \|f\| _{L^p _{r+\gamma/p}}
 \ge K_T \,
 \|f\|^{1 + \frac\gamma{pr}} _{L^p _r}, \]
where $K_T = \big( \sup_{t\in [0,T]} \| f\|_{L^p}(t) \big)^{-\frac\gamma{rq}}$. But this
last quantity is finite (thanks to estimate~\eqref{ghgh}).
We thus obtain the following a priori differential inequality on $\|f\|^p _{L^p _r}$:
 \[ \frac{d}{dt} \|f\|^p _{L^p _r} \le 
 - \, K_T \,  \big(\|f\|^p _{L^p _r}\big)^{1 + \frac\gamma{pr}}
 + C \, \|f\|^p _{L^p _r} \]
Using a standard argument (first used by Nash for parabolic equations) 
of comparison with the Bernouilli differential equation 
 \[ y' = -K_T \, y^{1+ \frac\gamma{pr}} + C \, y ,\]
whose solutions can be computed explicitely, we see that for all $0<t\le T$, 
 \[ \|f\| _{L^p_r}(t) < +\infty ,\]
more precisely 
 \begin{equation}\label{eq:upbdW}
 \|f\| _{L^p _r}(t) \le 
         \left[ \frac{C}{K_T \, \left(1-e^{-\frac{C \gamma}{pr} t}\right)} \right]^{r/\gamma} .
 \end{equation}
This concludes the proof of Theorem~\ref{theo:main} for local in times bounds. It remains 
to study  more accurately the behavior of these bounds when $t$ goes to infinity.
\bigskip

\Remarks
\smallskip

 1. Notice that the upper bound~\eqref{eq:upbdW} cannot be optimal since 
for example if  $\|f_0\|_{L^p _q} < +\infty$ then $\|f\|_{L^p _q} < +\infty$
uniformly on $[0,T]$ by the argument above, and the a priori differential 
inequality~\eqref{eq:apriori} implies that the quantity 
$\|f\|_{L^p _{q+\gamma/p}}$ is integrable at $t \sim 0^+$, which is not necessarily the case 
of the right-hand side term in~\eqref{eq:upbdW}. 
\smallskip

2. Note that in the previous computation, one should use approximate solutions 
of the Boltzmann equation in order to give a completely rigorous proof.
For example, solutions of the equation
 \begin{equation*}
 \left\{
 \begin{array}{ll} 
 \ \ \pa_t f_{\var} &= \ Q(f_{\var},f_{\var}) + \var\, \Delta_v f_{\var},\\
 f_{\var}(0,\cdot) &= \ f_{in} * \phi_{\var},
 \end{array}
 \right.
 \end{equation*}
where $\phi_{\var}$ is a sequence of mollifiers, can be used.
This point does not lead to any difficulties.
\smallskip 

3. It is also possible to get a slightly less stringent condition on the $L^1$ 
moments of the initial data $f_0$ by using the appearance of the $L^1$ moments of $f$ (in 
the case $\gamma > 0$). 
\medskip 

\section{Behavior for large times}\label{sec:unif}

The goal of this section is to conclude the proof of Theorem~\ref{theo:main} by showing
that the bounds on the $L^p$ moments  are uniform when $t\to +\infty$.

Our starting point is a stronger result than Proposition~\ref{prop:estim3}, which is 
a particular case of a result proven in~\cite{mv} (where the result holds 
for every collision kernel which satisfies angular integrability),  and is based on the 
regularity property of the gain term of the cutoff collision kernel. 
 This result writes:
\begin{proposition}[cf.~\cite{mv}, Theorem~4.1]\label{prop:estim4} 
 Let $B$ satisfy Assumptions~\eqref{eq:hyp1}, \eqref{eq:hyp2}, \eqref{eq:hyp3}. 
 We suppose moreover that 
 $b$ has its support in $[\theta_0,\pi/2]$. Then, 
 for all $p>1$, $q \ge 0$ and $f$ nonnegative with bounded entropy 
 and $L^1 _{2q + 2}$ norm, we have 
  \begin{equation}
  \int_{\R^N} Q(f,f)(v) \, f^{p-1}(v) \, \langle v \rangle ^{pq} \, dv \le 
  C^+ (b)  \, \|f\|^{p(1-\var)} _{L^p _q} 
  - K^- (b) \, \|f\|^p _{L^p _{q+\gamma/p}}
  \end{equation}
 with 
  \[ C^+ (b) = C^+ \, \left( \int_{\ens{S}^{N-1}} b \, d\sigma \right), \ \ \ \ 
     K^- (b) = K^- \, \left( \int_{\ens{S}^{N-1}} b \, d\sigma \right), \]
 and $C^+$, $K^-$ are positive constants. Both depend on an upper bound on 
 the entropy and the $L^1 _{2q + 2}$ norm of $f$ 
 and a lower bound on $\|f\|_{L^1}$; $C^+$ also depends on $\theta_0$.  
 Finally $\var \in (0,1)$ is a constant depending only on the 
 dimension $N$ and $p$.
 \end{proposition}
Gathering now Corollary~\ref{coro:estim2} with Proposition~\ref{prop:estim4}, we get  the
 \begin{proposition}\label{prop:estim6} 
 Let $B$ satisfy Assumptions~\eqref{eq:hyp1}, \eqref{eq:hyp2}, \eqref{eq:hyp3}, 
$p$ belong to $]1,+\infty[$ and $q\ge 0$. We suppose moreover that
 $pq \ge 2$ if $\nu \in (-2,-1]$ and $pq \ge 4$
 if $\nu \in (-3,-2]$.
 Then, for $f$ nonnegative with bounded entropy and $L^1 _{\max\{pq,2q\}+2}$ norm, we have
  \begin{equation}
  \int_{\R^N} Q(f,f)(v) \, f^{p-1}(v) \, \langle v \rangle ^{pq} \, dv \le 
  C^+ \, \|f\|^{p(1-\var)} _{L^p _q} 
  - K^- \, \|f\|^p _{L^p _{q+\gamma/p}}
  \end{equation}
 for some positive constants $C^+$ and $K^-$ depending on an upper bound on 
 $\|f\|_{L^1 _{\max\{pq,2q\}+2}}$, an upper bound on the entropy 
 and a lower bound on $\|f\|_{L^1}$. Finally $\var \in (0,1)$ is a constant depending only on the 
 dimension $N$ and $p$.
 \end{proposition}
\begin{proof}[Proof of Proposition~\ref{prop:estim6}]
The proof is exactly the same as that of Proposition~\ref{prop:estim5}. It is based 
on the splitting 
 \begin{equation*}
 b = b_c ^{\theta_0} + b_r ^{\theta_0}
 \end{equation*}
and the use of Corollary~\ref{coro:estim2} for
 \[ \int_{\R^N} Q_r (f,f)(v) \, f^{p-1}(v) \, \langle v \rangle ^{pq} \, dv \]
and Proposition~\ref{prop:estim4} for
 \[ \int_{\R^N} Q_c (f,f)(v) \, f^{p-1}(v) \, \langle v \rangle ^{pq} \, dv. \]
\end{proof} 

We now can prove that the bound on the $L^p$ moments is uniform for large times. 
Indeed, Proposition~\ref{prop:estim6} leads to the 
following a priori differential inequality on $y(t) = \|f(t,\cdot)\|^p _{L^p _q}$:
 \[ y' \le C \, y^{1-\var} - K \, y. \]
Then, by  a maximum principle, we see that $y(t)$ is bounded 
on $[\tau, +\infty[$ as soon as it is finite at time $\tau$.
The  explicit estimate is in fact:
 \[ \forall t\ge \tau, \ \ \ y(t) \le 
 \max \left\{ y(\tau); \, \left( \frac{C}{K} \right)^{1/\var} \right\}. \]
\bigskip 

\noindent
{\bf{Acknowledgment}}: Support by the European network HYKE, funded by the EC as
contract HPRN-CT-2002-00282, is acknowledged. 

\begin{flushleft} \signld \end{flushleft}
\vspace*{-55mm} \begin{flushright} \signcm \end{flushright}

\end{document}